\documentclass[letterpaper,12pt]{article}
\usepackage{amsmath,amssymb,amsthm}
\usepackage{graphicx}
\usepackage{hyperref}
\usepackage{subfig}
\usepackage{listings}

\newtheorem{theorem}{Theorem}
\newtheorem{example}{Example}
\newtheorem{lemma}[theorem]{Lemma}

\title{A construction of cospectral graphs\\ for the normalized Laplacian}
\author{Steve Butler\thanks{ This work was done with support of an NSF Mathematical Sciences Postdoctoral Fellowship while at UCLA.}\\
\small Department of Mathematics\\[-0.8ex]
\small Iowa State University\\[-0.8ex]
\small Ames, IA 50011\\[-0.8ex]
\small \texttt{butler@iastate.edu}
\and 
Jason Grout\\
\small Department of Mathematics and Computer Science\\[-0.8ex]
\small Drake University\\[-0.8ex]
\small Des Moines, IA  50311\\[-0.8ex]
\small \texttt {\tt jason.grout@drake.edu}}

\date{\small Mathematics Subject Classification: 05C50}

\begin{document}
\maketitle

\begin{abstract}
We give a method to construct cospectral graphs for the normalized Laplacian by a local modification in some graphs with special structure.  Namely, under some simple assumptions, we can replace a small bipartite graph with a cospectral mate without changing the spectrum of the entire graph.  We also consider a related result for swapping out biregular bipartite graphs for the matrix $A+tD$.

We produce (exponentially) large families of non-bipartite, non-regular graphs which are mutually cospectral, and also give an example of a graph which is cospectral with its complement but is not self-complementary.
\end{abstract}

\noindent{\bf Keywords:}~~normalized Laplacian; cospectral; bipartite subgraph swapping

\bigskip

\section{Introduction}

Spectral graph theory examines relationships between the structure of a graph and the eigenvalues (or spectrum) of a matrix associated with that graph.  Different matrices are able to give different information, but all the common matrices have limitations.  This is because there are graphs which have the same spectrum for a certain matrix but different structure---such graphs are called cospectral with respect to that matrix. 

\begin{figure}[t]
\centering
\hfil \subfloat[Cospectral for $A$]{\includegraphics[scale=1.2]{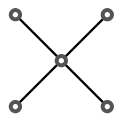}\hspace{.205in}\includegraphics[scale=1]{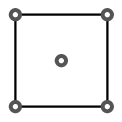}\label{fig:1a}}
\hfil\hfil \subfloat[Cospectral for $L$]{\includegraphics[scale=1.2]{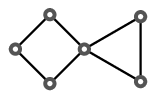}\hspace{.205in}\includegraphics[scale=1]{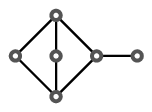}\label{fig:1b}}
\hfil\hfil \subfloat[Cospectral for $Q$]{\includegraphics[scale=1.2]{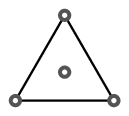}\hspace{.205in}\includegraphics[scale=1]{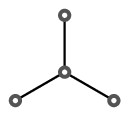}\label{fig:1c}}
\hfil\hfil \subfloat[Cospectral for $\mathcal{L}$]{\includegraphics[scale=1.2]{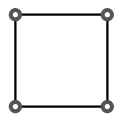}\hspace{.205in}\includegraphics[scale=1]{SpecPics-6}\label{fig:1d}} \hfil \hfil
\caption{Examples of cospectral graphs for $A$, $L$, $Q$ and $\mathcal{L}$.}
\label{fig:co}
\end{figure}

The following are some of the matrices studied in spectral graph theory:
\begin{itemize}
\item The adjacency matrix $A$.  This is defined by $A(u,v)=1$ when $u$ and $v$ are adjacent and $0$ otherwise.  The spectrum of the adjacency matrix can determine the number of edges and if a graph is bipartite, but it cannot determine if a graph is connected (see \cite{CDS} for more information).  See Figure~\ref{fig:1a}, the ``Saltire pair'' \cite{Cvetkovic,HaemersSpence}, for an example of two graphs which are cospectral with respect to the adjacency matrix but not any of the other matrices we will discuss.
\item The Laplacian $L=D-A$.  This matrix is also known as the combinatorial Laplacian and is found by taking the difference of the diagonal degree matrix $D$ and the adjacency matrix.  The spectrum of the Laplacian can determine the number of edges and the number of connected components, but it cannot determine if a graph is bipartite (see \cite{CDS} for more information).  See Figure~\ref{fig:1b}  \cite{Dam,HaemersSpence}, for an example of two graphs which are cospectral with respect to the Laplacian but not any of the other matrices we will discuss.
\item The signless Laplacian $Q=D+A$.  This matrix is found by taking the sum of the diagonal degree matrix $D$ and the adjacency matrix.  The spectrum of the signless Laplacian can determine the number of edges and the number of connected components which are bipartite, but it cannot determine if a graph is bipartite or connected (see \cite{sign1,sign2,sign3} for more information).  See Figure~\ref{fig:1c}  \cite{HaemersSpence}, for an example of two graphs which are cospectral with respect to the signless Laplacian but not any of the other matrices we will discuss.
\item The normalized Laplacian ${\cal L}$.  This matrix is defined by ${\cal L}=D^{-1/2}LD^{-1/2}$, where by convention if we have an isolated vertex then it will contribute $0$ to the spectrum.  The spectrum of the normalized Laplacian is closely related to the spectrum of the probability transition matrix of a random walk.  This spectrum can determine if a graph is bipartite and the number of connected components, but it cannot determine the number of edges (see \cite{Chung} for more information).  See Figure~\ref{fig:1d} for an example of two graphs which are cospectral with respect to the normalized Laplacian but not any of the other matrices we will discuss.
\end{itemize}

One way to understand what structure the spectrum of a matrix cannot identify is to study cospectral graphs.  Cospectral graphs for the adjacency matrix \cite{DamHaemersKoolen,GodsilMcKay,GodsilMcKay2,HaemersSpence,JohnsonNewman} and the Laplacian matrix \cite{HaemersSpence,Merriss,Tan} have been studied, particularly for graphs with few vertices.  Cospectral graphs for the signless Laplacian have been little studied beyond their enumeration and rules which apply to all of these matrices \cite{HaemersSpence}.

Little is also known about cospectral graphs for the normalized Laplacian compared to other matrices.  Previously, the only cospectral graphs were bipartite (complete bipartite graphs \cite{Chung,Tan} and bipartite graphs found by ``unfolding'' a small bipartite graph in two ways \cite{Butler}) or were regular and cospectral for $A$ (since a regular graph is cospectral for all of these four matrices if it is cospectral for any one of them). 

In Table~\ref{table}, we have listed the number of graphs with a cospectral mate for the normalized Laplacian matrix for graphs on nine or fewer vertices, counted using Sage (see Appendix~\ref{sagecode} for example code).  The number of cospectral graphs for the normalized Laplacian have also been counted by Wilson and Zhu \cite{WilsonZhu}, though they give percentages and not the count of how many have cospectral mates.  We have also included similar counts for the other three matrices which come from \cite{brouwer,HaemersSpence}.

\begin{table}[ht]
\[
\begin{array}{r|r|r|r|r|r}
\mbox{\#vertices}&\mbox{\#graphs}&A&L{=}D{-}A&Q{=}D{+}A&\mathcal{L}\\ \hline
 1&         1&        0&       0&       0&    0\\
 2&         2&        0&       0&       0&    0\\
 3&         4&        0&       0&       0&    0\\
 4&        11&        0&       0&       2&    2\\
 5&        34&        2&       0&       4&    4\\
 6&       156&       10&       4&      16&   14\\
 7&      1044&      110&     130&     102&   52\\
 8&     12346&     1722&    1767&    1201&  201\\
 9&    274668&    51039&   42595&   19001& 1092
\end{array}
\]
\caption{Number of graphs with a cospectral mate for the various matrices}
\label{table}
\end{table}

\begin{figure}[hb]
\centering
\includegraphics[scale=1]{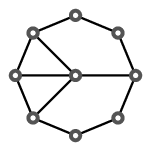}\hspace{0.5in}\includegraphics[scale=1]{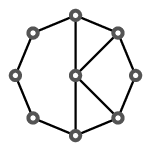}
\caption{An example of two cospectral graphs for the adjacency matrix related by switching.}
\label{fig:switch}
\end{figure}

Given the relatively small number of graphs with a cospectral mate with respect to the normalized Laplacian, it is surprising that so little is known about forming cospectral graphs for that matrix.  The problem is that some of the main tools that are used to form cospectral graphs for other matrices do not generalize to the normalized Laplacian.  One such example is a technique
known as switching, which is accomplished by replacing edges by non-edges and non-edges by edges between two subsets (see \cite{GodsilMcKay2, HaemersSpence, Seidel}; a simple example is shown in Figure~\ref{fig:switch}).  Given some basic assumptions, this is an easy way to construct cospectral graphs for the adjacency matrix.  However, switching does not in general work for the normalized Laplacian, in particular, it will only be guaranteed to work when the degrees are unchanged (see \cite{Cavers}).  So the graphs shown in Figure~\ref{fig:switch} are {\em not}\/ cospectral with respect to the normalized Laplacian.

In Section~\ref{sec:swap} we will introduce a method to construct cospectral graphs for the normalized Laplacian.  This construction will work similarly to switching in that we will make a small local change to the graph by swapping in one bipartite graph with a cospectral mate and show the two graphs still share the same eigenvalues.  In Section~\ref{sec:other} we will show that if we add additional constraints to the bipartite graphs which are swapped then the resulting graphs are also cospectral with respect to $A$, $L$ and $Q$.  In Section~\ref{sec:applications} we will show how to use this construction to produce large families of graphs which are mutually cospectral.  Finally, in Section~\ref{sec:comments} we will give some concluding remarks.

\section{Swapping bipartite subgraphs}\label{sec:swap}
The method of finding cospectral graphs for the adjacency matrix by switching reduces to making a local change of the graph from one bipartite subgraph to its complement and showing that the spectrum is unchanged.  We will consider something similar, namely a local change, but instead of replacing a bipartite subgraph with its complement we will swap out a bipartite subgraph for a cospectral mate.

For a subset $W$ of the vertices $V$ of $G$, we will let $G[W]$ be the induced subgraph of $G$ on the vertex set $W$.

\begin{theorem}\label{thm:biswitch}
Let $P_1$ and $P_2$ be bipartite, cospectral graphs with respect to the normalized Laplacian on the vertex set $B\cup C$ such that all edges go between $B$ and $C$ and where all vertices in $B$ have degree $k$ for both graphs.

Let $G_1$ be a graph on the vertices $A\cup A'\cup B\cup C$ where $G_1[A\cup A']$ is an arbitrary graph; $G_1[B]$ and $G_1[C]$ have no edges; there are no edges going between $A$ and $B$, between $A$ and $C$, or between $A'$ and $C$; $G_1[A'\cup B]$ is a complete bipartite graph; and $G_1[B\cup C]=P_1$.  The graph $G_2$ is defined similarly except that $G_2[B\cup C]=P_2$.  Then $G_1$ and $G_2$ are cospectral with respect to the normalized Laplacian.

If the dimension of the eigenspace associated with $\lambda=1$ intersected with the subspace of vectors that are nonzero only on $B$ is the same for both $P_1$ and $P_2$, then the graphs $H_1$ and $H_2$ are also cospectral with respect to the normalized Laplacian, where $H_1$ only differs from $G_1$ in that $H_1[B]$ is the complete graph and similarly  $H_2$ only differs from $G_2$ in that $H_2[B]$ is the complete graph.
\end{theorem}

\begin{figure}[hftb]
\centering
\subfloat[$G_1$]{\includegraphics{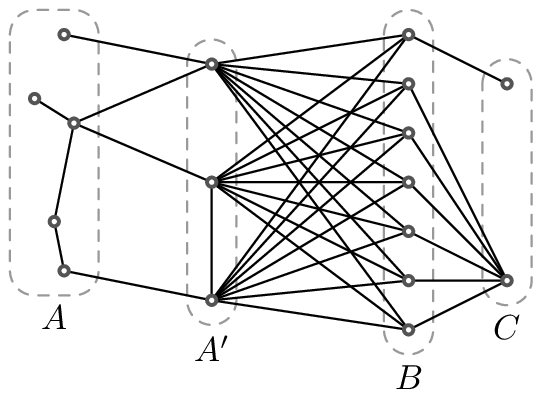}}\hfil
\subfloat[$G_2$]{\includegraphics{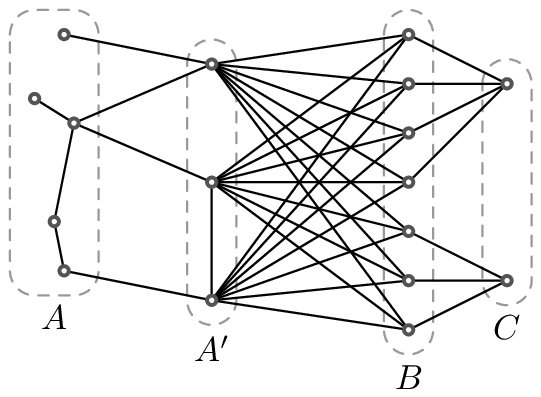}}
\caption{An example of cospectral graphs using the construction given in Theorem~\ref{thm:biswitch}.}
\label{fig:exswitch}
\end{figure}

An example of the construction described in Theorem~\ref{thm:biswitch} is shown in Figure~\ref{fig:exswitch}. 
In this case $P_1$ and $P_2$ are the graphs $K_{1,1}\cup K_{1,6}$ and $K_{1,4}\cup K_{1,3}$ where we have put the degree $1$ vertices into $B$.  It is well known that the spectrum of a complete bipartite graph $K_{p,q}$ is $0^{[1]},1^{[p+q-2]},2^{[1]}$, where the exponent indicates multiplicity, hence these graphs are easily seen to be cospectral.  Moreover it is easy to check that the dimension of the eigenspace for the eigenvalue $1$ intersected with the subspace of vectors which are nonzero only on $B$ is $5$ for both graphs, i.e., we could also have put a complete graph on the vertices of $B$ in the graphs shown in Figure~\ref{fig:exswitch} and still had a cospectral pair.

\begin{figure}[hftb]
\centering
\subfloat[$P_1$]{\includegraphics{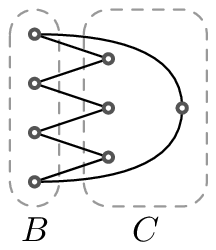}}\hfil
\subfloat[$P_2$]{\includegraphics{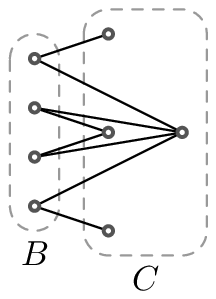}}
\caption{A pair of cospectral graphs satisfying the conditions of Theorem~\ref{thm:biswitch}.}
\label{fig:other1}
\end{figure}

\begin{figure}[hftb]
\centering
\subfloat[$P_1$]{\includegraphics{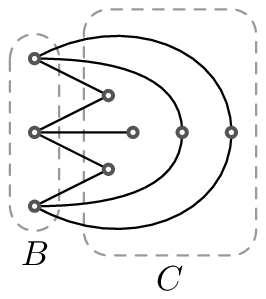}}\hfil
\subfloat[$P_2$]{\includegraphics{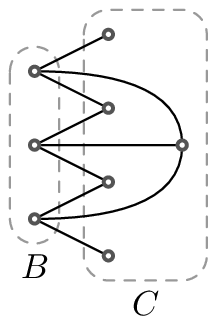}}
\caption{A pair of cospectral graphs satisfying the conditions of Theorem~\ref{thm:biswitch}.}
\label{fig:other2}
\end{figure}

Additional examples of graphs $P_1$ and $P_2$ satisfying (both) the conditions on Theorem~\ref{thm:biswitch} are $P_1=K_{\ell,p}\cup K_{\ell,q}$ and $P_2=K_{\ell,r}\cup K_{\ell,s}$ where the vertices of degree $\ell$ are all placed in $B$ and $p+q=r+s$.  Two further examples are shown in Figures~\ref{fig:other1} and~\ref{fig:other2}.

To prove Theorem~\ref{thm:biswitch} we will find it convenient to work with the harmonic eigenvectors of the normalized Laplacian.  Namely, if ${\bf x}\neq{\bf 0}$ is an eigenvector associated with the eigenvalue $\lambda$, i.e., ${\cal L}{\bf x}=\lambda{\bf x}$, then the corresponding harmonic eigenvector is ${\bf y}=D^{-1/2}{\bf x}$.  This translates the relationship ${\cal L}{\bf x}=\lambda{\bf x}$ into $(D-A){\bf y}=\lambda D{\bf y}$, which at a vertex $v$ becomes
\begin{equation}\label{eq:harmonic}
\sum_{u{\sim}v}{\bf y}(u)=(1-\lambda){\bf y}(v)d(v),
\end{equation}
where $d(v)$ is the degree of the vertex $v$.  We will say that two harmonic eigenvectors ${\bf y}_1$ and ${\bf y}_2$ are orthogonal if ${\bf y}_2^*D{\bf y}_1=0$, i.e., if the corresponding eigenvectors are orthogonal.

Before we begin the proof of Theorem~\ref{thm:biswitch} it will be useful to make some observations about bipartite graphs that will come up in the proof.  

\begin{lemma}\label{lem:help}
Let $P$ be a bipartite graph on the vertices $B\cup C$ where all edges go between $B$ and $C$ and the vertices in $B$ all have degree $k$.  Further, if ${\bf x}$ is an eigenvector for the normalized Laplacian then we can write it as ${\bf x}={\bf b}+{\bf c}$ where ${\bf b}$ is the vector ${\bf x}$ restricted to the vertices of $B$, and similarly for ${\bf c}$.  Then the following hold:
\begin{itemize}
\item[(a)] If ${\cal L}{\bf x}={\bf x}$, i.e., ${\bf x}$ is an eigenvector for the eigenvalue $1$, then ${\cal L}{\bf b}={\bf b}$ and ${\cal L}{\bf c}={\bf c}$.  In other words, we can divide the eigenspace for the eigenvalue $1$ between those which are nonzero only on the vertices of $B$ and nonzero only on the vertices of $C$.
\item[(b)] If ${\bf x}={\bf b}+{\bf c}$ is an eigenvector associated with $\lambda$, then ${\bf b}-{\bf c}$ is an eigenvector associated with $2-\lambda$.
\item[(c)] If ${\bf x}_1={\bf 1}_B+{\bf 1}_C,\ldots,{\bf x}_i={\bf b}_i+{\bf c}_i,\ldots,{\bf x}_m={\bf 1}_B-{\bf 1}_C$ are a set of orthogonal harmonic eigenvectors of $P$, then for $1<i<m$ we have ${\bf b}_i$ is orthogonal to ${\bf 1}_B$ and ${\bf c}_i$ is orthogonal to ${\bf 1}_C$.  In particular,  for $1<i<m$ the sum of the entries of ${\bf b}_i$ is $0$.
\end{itemize}
\end{lemma}

Parts (a) and (b) easily follow from \eqref{eq:harmonic}.  For (c) we note that ${\bf x}_i$ will be orthogonal to the vectors ${1\over2}({\bf x}_1+{\bf x}_m)={\bf 1}_B$ and ${1\over2}({\bf x}_1-{\bf x}_m)={\bf 1}_C$, also since the vertices in $B$ are regular the result on the sum of the entries of ${\bf b}_i$ follows by the definition of orthogonality.

\begin{proof}[Proof of Theorem~\ref{thm:biswitch}]
Let ${\bf x}_1={\bf 1}_B+{\bf 1}_C,\ldots,{\bf x}_i={\bf b}_i+{\bf c}_i,\ldots,{\bf x}_m={\bf 1}_B-{\bf 1}_C$ be a full set of orthogonal harmonic eigenvectors of $P_1$ associated with $\lambda_1=0,\ldots,\lambda_i,\ldots,\lambda_m=2$ respectively.  Let $s$ denote the degree of a vertex in $B$ for the graph $G_1$ and $k$ the degree in $P_1$.  Then for $1<i<m$ we have that $\widehat{\bf x}_i={\bf b}_i+\sqrt{s\over k}{\bf c}_i$ is a harmonic eigenvector for $G_1$ associated with eigenvalue $\gamma_i=1-(1-\lambda_i)\sqrt{k\over s}$.  To see this we need to examine what happens for each vertex $v$ in $G_1$ using \eqref{eq:harmonic}.
\begin{itemize}
\item For $v\in A$: Then all of the entries of the vertex and its neighbors are all $0$ and so \eqref{eq:harmonic} trivially holds.
\item For $v\in A'$: All of the neighbors in $A$ and $A'$ are $0$ while the sum of the entries in $B$ will be $0$ by Lemma~\ref{lem:help}(c), so both sides again are $0$ and \eqref{eq:harmonic} holds.
\item For $v\in B$: We note that all of the nonzero elements of $\widehat{\bf x}_i$ adjacent to $v$ are in $C$, and further we have that for the harmonic eigenvector ${\bf x}_i$ for $P_1$ that
\[
\sum_{u{\sim}v}{\bf x}_i(u)=(1-\lambda_i){\bf x}_i(v)k.
\]
So we have for $G_1$ that 
\begin{multline*}
\sum_{u{\sim}v}\widehat{\bf x}_i(u)=\sqrt{s\over k}\sum_{u{\sim}v}{\bf x}_i(u)=\sqrt{s\over k}(1-\lambda_i){\bf x}_i(v)k=\sqrt{k\over s}(1-\lambda_i)\widehat{\bf x}_i(v)d(v)\\=(1-\gamma_i)\widehat{\bf x}_i(v)d(v),
\end{multline*}
so \eqref{eq:harmonic} again holds.
\item For $v\in C$: We proceed similarly as we did for vertices in $B$ and we get
\[
\sum_{u{\sim}v}\widehat{\bf x}_i(u)=\sum_{u{\sim}v}{\bf x}_i(u)=(1-\lambda_i){\bf x}_i(v)d(v)=\sqrt{k\over s}(1-\lambda_i)\widehat{\bf x}_i(v)d(v)=(1-\gamma_i)\widehat{\bf x}_i(v)d(v).
\]
\end{itemize}

Next we note by Lemma~\ref{lem:help}(c) that ${\bf 1}_B$ and ${\bf 1}_C$ are orthogonal to each of the vectors $\widehat{\bf x}_i$ for $1<i<m$.  On $P_1$, the dimension of the space orthogonal to the $\widehat{\bf x}_i$ for $1<i<m$ restricted to $P_1$ is two.  A basis for the orthogonal complement of $\mathrm{span}\{\widehat{\bf x}_i\}$ is $\{{\bf 1}_B,{\bf 1}_C\}$, so any harmonic vector orthogonal to all of these {\em must}\/ be a linear combination of ${\bf 1}_B$ and ${\bf 1}_C$.  In particular, if we let ${\bf y}$ be a harmonic eigenvector of $G_1$ that is orthogonal to all of the $\widehat{\bf x}_i$ for $1<i<m$, then ${\bf y}$ restricted to $B$ and $C$ is orthogonal to the $\widehat{\bf x}_i$, so we have ${\bf y}|_{B\cup C}=b{\bf 1}_B+c{\bf 1}_C$ for some constants $b$ and $c$.  In other words, the harmonic eigenvector ${\bf y}$ is constant on the vertices of $B$ and constant on the vertices of $C$.

Everything that we have done for $G_1$ carries over to $G_2$.  In particular, since $P_1$ and $P_2$ are cospectral, then the $\gamma_i$ found by generalizing the harmonic eigenvectors of $P_1$ and $P_2$ remain the same.  Furthermore, every other harmonic eigenvector orthogonal to the $\widehat{\bf x}_i$ must be constant on the vertices of $B$ and the vertices of $C$.  To finish off the first claim of the theorem, we now only need to observe that {\em any}\/ harmonic eigenvector which is orthogonal to all of the $\widehat{\bf x}_i$ in $G_1$ is also a harmonic eigenvector for $G_2$ for the same eigenvalue.  Again to see this we only need to consider what happens for each vertex.
\begin{itemize}
\item For $v\in A\cup A'$: Then all of the entries of the vertex and its neighbors are the same for both graphs, so \eqref{eq:harmonic} trivially holds.
\item For $v\in B$: All of the neighbors in $A'$ are the same, and we have the exact same number of neighbors in $C$ as before with the same value and so again \eqref{eq:harmonic} holds.
\item For $v\in C$: For the graph $G_1$, let $\beta$ be the fixed value of the vertices in $B$.  We have
\[
\sum_{u{\sim}v}{\bf y}(u)=d(v)\beta=(1-\lambda){\bf y}(v)d(v).
\]
Of course the $d(v)$ terms cancel and we are left with $\beta=(1-\lambda){\bf y}(v)$.  So even though the degree of the vertex in $C$ might change, it will have no effect on this relationship and so \eqref{eq:harmonic} again holds.
\end{itemize}

In summary, we were able to find $m-2$ harmonic eigenvectors for $G_1$ and $G_2$ that gave the same set of eigenvalues.  For any other harmonic eigenvector orthogonal to these, the same harmonic eigenvector worked for both graphs and so the remaining set of eigenvalues also agreed.  So we can conclude that $G_1$ and $G_2$ are cospectral with respect to the normalized Laplacian.

\bigskip

We now turn to the second statement of the theorem.  Let $s$ be the degree of a vertex in $B$ in the graph $H_1$.  So let ${\bf x}_1={\bf 1}_B+{\bf 1}_C,\ldots,{\bf x}_i={\bf b}_i+{\bf c}_i,\ldots,{\bf x}_m={\bf 1}_B-{\bf 1}_C$ be a full set of orthogonal harmonic eigenvectors of $P_1$ associated with $\lambda_1=0,\ldots,\lambda_i,\ldots,\lambda_m=2$ respectively. First, we consider the harmonic eigenvectors of the graph $P_1$ in the eigenspace corresponding to $\lambda=1$.  By Lemma~\ref{lem:help}(a) we can assume that the vectors ${\bf x}_i$ are either of the form ${\bf b}_i$ or ${\bf c}_i$.  Let the corresponding harmonic eigenvector of $H_1$ be $\widehat{\bf x}_i={\bf x}_i$, i.e., we simply expand ${\bf x}_i$ to be zero outside of $P_1$.  Then we have two cases.
\begin{itemize}
\item If ${\bf x}_i={\bf b}_i$:  In this case, we know there are no problems for the vertices in $A$ (every term in \eqref{eq:harmonic} is $0$) or $C$ (since \eqref{eq:harmonic} reduces to what was done in $P_1$).  For the vertices in $A'$ we can use Lemma~\ref{lem:help}(c) to see that both sides of \eqref{eq:harmonic} are $0$.  Finally, suppose that $v$ is a vertex in $B$.  Then we have
\[
\sum_{u{\sim}v}\widehat{\bf x}_i(u)=\sum_{u\in B}\widehat{\bf x}_i(u)-\widehat{\bf x}_i(v)=-\widehat{\bf x}_i(v)=-{1\over s}\widehat{\bf x}_i(v)s=-{1\over s}\widehat{\bf x}_i(v)d(v).
\]
In particular, we have that $1-\gamma_i=-{1\over s}$ is an eigenvalue, i.e., that $\gamma_i=1+{1\over s}$ is an eigenvalue, for this harmonic eigenvector.
\item If ${\bf x}_i={\bf c}_i$:  In this case we similarly know there are no problems for vertices in $A$, $A'$ and $B$.  For vertices in $C$, \eqref{eq:harmonic} reduces to what we had in $P_1$ and so we can conclude that this is a harmonic eigenvector for the eigenvalue $\gamma_i=1$.
\end{itemize}

Now suppose that ${\bf x}_i={\bf b}_i+{\bf c}_i$ is associated with $\lambda_i\neq 1$ for $P_1$.  Then we now will create two harmonic eigenvectors for $H_1$, namely
\begin{align*}
{\bf y}_i^1&={\bf b}_i+t_1{\bf c}_i~=~{\bf b}_i+\bigg({1+\sqrt{1+4(1-\lambda_i)^2sk}\over 2k(1-\lambda_i)}\bigg){\bf c}_i,\\
{\bf y}_i^2&={\bf b}_i+t_2{\bf c}_i~=~{\bf b}_i+\bigg({1-\sqrt{1+4(1-\lambda_i)^2sk}\over 2k(1-\lambda_i)}\bigg){\bf c}_i.
\end{align*}
These will be associated with the eigenvalues of
\[
\gamma_1={2s+1-\sqrt{1+4(1-\lambda_i)^2sk}\over 2s}\qquad\mbox{and}\qquad
\gamma_2={2s+1+\sqrt{1+4(1-\lambda_i)^2sk}\over 2s},
\]
respectively, for $H_1$.  Some simple computations show that the following relationships hold:
\[
t_1(1-\gamma_1)=t_2(1-\gamma_2)=1-\lambda_i,~~
t_1(1-\lambda_i)k-1=(1-\gamma_1)s,~~\mbox{and}~~
t_2(1-\lambda_i)k-1=(1-\gamma_2)s.
\]

(Before we proceed to the next step of showing that these are indeed harmonic eigenvectors and eigenvalues we first make the observation that we are not creating more harmonic eigenvectors than we had before.  This is because in Lemma~\ref{lem:help}(b) we can pair up harmonic eigenvectors and eigenvalues.  In particular, we would have generated the same new harmonic eigenvectors and eigenvalues if we had used $2-\lambda_i$ and ${\bf b}_i-{\bf c}_i$.  So really we are taking pairs of eigenvectors and eigenvalues to new pairs of eigenvectors and eigenvalues.)

To verify that these are harmonic eigenvectors with the specified eigenvalues we need to examine what happens for each vertex $v$ in $H_1$ using \eqref{eq:harmonic}.  We will step through ${\bf y}_i^1$ and $\gamma_1$, the arguments for ${\bf y}_i^2$ and $\gamma_2$ are the same.
\begin{itemize}
\item For $v\in A$: Then all of the entries of the vertex and its neighbors are all $0$ and so \eqref{eq:harmonic} trivially holds.
\item For $v\in A'$: All of the neighbors in $A$ and $A'$ are $0$ while the sum of the entries in $B$ will be $0$ by Lemma~\ref{lem:help}(c), so both sides again are $0$ and \eqref{eq:harmonic} holds.
\item For $v\in B$: We note that all of the nonzero elements of $\widehat{\bf x}_i$ adjacent to $v$ are in $B$ and $C$, and further we have that for the harmonic eigenvector ${\bf x}_i$ for $P_1$ that
\[
\sum_{u{\sim}v\atop u\in C}{\bf x}_i(u)=(1-\lambda_i){\bf x}_i(v)k.
\]
So we have for $H_1$ that 
\begin{multline*}
\sum_{u{\sim}v}{\bf y}_i^1(u)=t_1\sum_{u{\sim}v\atop u\in C}{\bf x}_i(u)+\sum_{u\in B}{\bf x}_i(u)-{\bf x}_i(v)=t_1(1-\lambda_i){\bf x}_i(v)k-{\bf x}_i(v)\\
=\big(t_1(1-\lambda_i)k-1\big){\bf x}_i(v)=(1-\gamma_1){\bf y}_i^1(v)s=(1-\gamma_1){\bf y}_i^1(v)d(v).
\end{multline*}
so \eqref{eq:harmonic} again holds.
\item For $v\in C$: We proceed similarly as we did for vertices in $B$ and we get
\[
\sum_{u{\sim}v}{\bf y}_i^1(u)=\sum_{u{\sim}v}{\bf x}_i(u)=(1-\lambda_i){\bf x}_i(v)d(v)=t_1(1-\gamma_1){\bf x}_i(v)d(v)=(1-\gamma_1){\bf y}_i^1(v)d(v).
\]
\end{itemize}

The remainder of this case now proceeds as before.  Namely, everything that we did for $H_1$ carries over for $H_2$ (counting multiplicity of eigenvalues).  Further, any other harmonic eigenvector orthogonal to the ones given must be constant on $B$ and $C$ and so any other harmonic eigenvector which works for $H_1$ also works for $H_2$.  Therefore we can conclude that they have the same set of eigenvalues, i.e., are cospectral with respect to the normalized Laplacian.
\end{proof}

\section{Swapping biregular bipartite subgraphs}\label{sec:other}
The graphs generated by Theorem~\ref{thm:biswitch} will generally not give cospectral graphs for other matrices.  So it is instructive to examine the proof and try to understand the point at which the fact that we were using the normalized Laplacian as compared to some other matrix came into play.  The key is understanding how \eqref{eq:harmonic} remains true even when the degrees of vertices in $C$ change---if the degree of a vertex in $C$ changes, the $d(v)$ term on the right side of \eqref{eq:harmonic} will change proportionally to the sum on the left so that the equality still holds.

In the case when the degrees of vertices do not change, for example, when $P_1$ and $P_2$ are both $(k,\ell)$-biregular, then the proof generalizes.  Recall that a bipartite graph is $(k,\ell)$-biregular if the vertices can be partitioned into $B\cup C$ where all edges go between $B$ and $C$, the vertices in $B$ all have degree $k$, and the vertices in $C$ have degree $\ell$.  We now have the following theorem.

\begin{theorem}\label{thm:biregswitch}
Let $P_1$ and $P_2$ be bipartite, cospectral with respect to the adjacency matrix, $(k,\ell)$-biregular graphs on the vertex set $B\cup C$ such that all vertices in $B$ have degree $k$ and all vertices in $C$ have degree $\ell$ and edges go between $B$ and $C$.  Further, the dimension of the eigenspace associated with $\lambda=0$ for the adjacency matrix intersected with the subspace of vectors which are nonzero only on $B$ is the same for both $P_1$ and $P_2$.  

Let $G_1$ be a graph on the vertices $A\cup A'\cup B\cup C$ where $G_1[A\cup A']$ is an arbitrary graph; $G_1[B]$ and $G_1[C]$ have no edges; there are no edges going between $A$ and $B$, between $A$ and $C$, or between $A'$ and $C$; $G_1[A'\cup B]$ is a complete bipartite graph; and $G_1[B\cup C]=P_1$.  The graph $G_2$ is defined similarly except that $G_2[B\cup C]=P_2$.  Then $G_1$ and $G_2$ are cospectral with respect to the matrix $A+tD$ for $t$ arbitrary.

Similarly, the graphs $H_1$ and $H_2$ are also cospectral with respect to $A+tD$ for $t$ arbitrary, where $H_1$ only differs from $G_1$ in that $H_1[B]$ is the complete graph and similarly  $H_2$ only differs from $G_2$ in that $H_2[B]$ is the complete graph.
\end{theorem}

Note when $t=0$ then $A+tD=A$, when $t=1$ then $A+tD=A+D=Q$, and when $t=-1$ then $A+tD=A-D=-L$, so that any such pairs of graphs are cospectral with respect to all of the matrices $A$, $L$, and $Q$.  Further, these are also cospectral for $\cal{L}$ since if $t=-\lambda-1$, we have
\begin{align*}
  \det (A+(-\lambda-1)D)&=\det((A-D)-\lambda D)\\
  &= \det (D^{1/2}(\mathcal{L}-\lambda I)D^{1/2})\\
  &= (\prod d_i)\det(\mathcal{L}-\lambda I).
\end{align*}
Since the two graphs have the same degree sequence, the product $\prod d_i$ is the same, so the characteristic polynomials also are the same.

In order to use this theorem, we must find two $(k,\ell)$-biregular bipartite graphs which are cospectral and for which the dimensions of the eigenspace of $0$ restricted to $B$ agree on both graphs.  There is one special case for which this is much easier, namely the theorem does not prohibit the possibility that $P_1$ and $P_2$ are the same graph (so trivially are cospectral so we only are reduced to checking the dimension of the eigenspaces).  This is not so trivial as it might sound since we are distinguishing the two parts of the bipartite graph in the proof of the theorem, so that while $P_1$ and $P_2$ are the same graph we are attaching them in two different ways and so the resulting graphs might not be isomorphic.

\begin{figure}[ht]
\centering
\subfloat[$P_1$]{\includegraphics{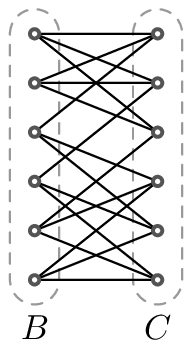}}\hfil
\subfloat[$P_2$]{\includegraphics{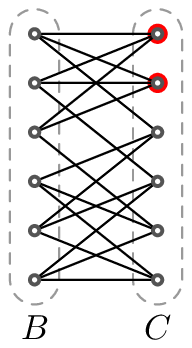}}
\caption{A pair of cospectral graphs satisfying the conditions of Theorem~\ref{thm:biregswitch}.}
\label{fig:2bi}
\end{figure}

As an example of this consider the graphs shown in Figure~\ref{fig:2bi}.  These are both the same graph (just flipped) and a simple check shows that the dimension of the eigenspace associated with $0$ and restricted to $B$ is $2$ in both graphs.  Therefore Theorem~\ref{thm:biregswitch} applies, and any two graphs to which these two graphs are attached are cospectral.  But note that in $P_2$ there are two vertices in $C$ with the same neighbors (marked), but there are no such pairs in $P_1$.  So the resulting graphs are {\em not} isomorphic as long as $A'\ne \emptyset$.

\begin{figure}[ht]
\centering
\subfloat[$H_1$]{\includegraphics{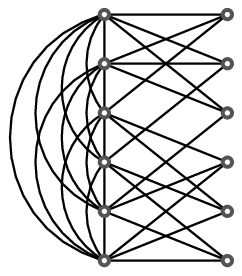}}\hfil
\subfloat[$H_2$]{\includegraphics{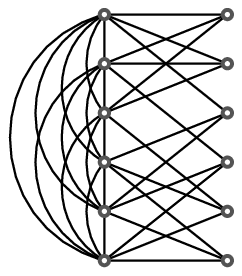}}
\caption{A pair of non-isomorphic cospectral graphs for $A$, $L$, $Q$ and ${\cal L}$.  (Note $\overline{H_1}=H_2$.)}
\label{fig:bicom}
\end{figure}

So for example if we consider the graphs $H_1$ and $H_2$ where $A\cup A'=\emptyset$ and we induce a complete graph on the vertices of $B$, then Theorem~\ref{thm:biregswitch} shows that the resulting graphs are cospectral, and again since there are two vertices with degree $3$ sharing common neighbors in one graph but not the other they are non-isomorphic.  These graphs are shown in Figure~\ref{fig:bicom}.  These graphs also happen to be complements of one another, and so give an example of a graph which is cospectral with its complement but not self-complementary.  (Note the graph shown in Figure~\ref{fig:1c} is another example of a graph which is cospectral with its complement, but in that case it is only cospectral with respect to $Q$ and not with respect to $A$, $L$ or ${\cal L}$.)

The proof for Theorem~\ref{thm:biregswitch} will proceed similar to the proof for Theorem~\ref{thm:biswitch}, so we will skip some of the routine computations and provide an outline of the proof.  Note that if ${\bf x}$ is an eigenvalue associated with eigenvalue $\lambda$, i.e., $(A+tD){\bf x}=\lambda{\bf x}$, then at a vertex $v$ we have
\begin{equation}\label{eq:t}
\sum_{u{\sim}v}{\bf x}(u)=\big(\lambda-td(v)\big){\bf x}(v).
\end{equation}
We first start with some simple properties of bipartite graphs that will be helpful in the proof.  The proofs are similar to Lemma~\ref{lem:help} and we will omit them.

\begin{lemma}\label{lem:bihelp}
Let $P$ be a bipartite $(k,\ell)$-biregular graph on the vertices $B\cup C$ where all edges go between $B$ and $C$.  Further, if ${\bf x}$ is an eigenvector for the adjacency matrix, then we can write it as ${\bf x}={\bf b}+{\bf c}$ where ${\bf b}$ is the vector ${\bf x}$ restricted to the vertices of $B$, and similarly for ${\bf c}$.  Then the following hold:
\begin{itemize}
\item[(a)] If $A{\bf x}={\bf 0}$, i.e., ${\bf x}$ is an eigenvector for the eigenvalue $0$, then $A{\bf b}={\bf 0}$ and $A{\bf c}={\bf 0}$.  In other words, we can divide the eigenspace for the eigenvalue $0$ between those which are nonzero only on the vertices of $B$ and nonzero only on the vertices of $C$.
\item[(b)] If ${\bf x}={\bf b}+{\bf c}$ is an eigenvector associated with $\lambda$, then ${\bf b}-{\bf c}$ is an eigenvector associated with $-\lambda$.
\item[(c)] If ${\bf x}_1=\sqrt k{\bf 1}_B+\sqrt \ell{\bf 1}_C,\ldots,{\bf x}_i={\bf b}_i+{\bf c}_i,\ldots,{\bf x}_m=\sqrt k{\bf 1}_B-\sqrt \ell{\bf 1}_C$ are a set of orthogonal eigenvectors of $P$, then for $1<i<m$ we have ${\bf b}_i$ is orthogonal to ${\bf 1}_B$ and ${\bf c}_i$ is orthogonal to ${\bf 1}_C$.  In particular,  for $1<i<m$ the sum of the entries of ${\bf b}_i$ is $0$.
\end{itemize}
\end{lemma}

\begin{proof}[Proof of Theorem~\ref{thm:biregswitch}]
We will first consider the case when the vertices on $B$ induce an empty graph.

Let ${\bf x}_1=\sqrt k{\bf 1}_B+\sqrt \ell{\bf 1}_C,\ldots,{\bf x}_i={\bf b}_i+{\bf c}_i,\ldots,{\bf x}_m=\sqrt k{\bf 1}_B-\sqrt \ell{\bf 1}_C$ be a full set of orthogonal eigenvectors of $P_1$ for the adjacency matrix associated with the eigenvalues $\lambda_1=\sqrt{k\ell},\ldots,\lambda_i,\ldots,-\sqrt{k\ell}$ respectively.  Let $s$ denote the degree of a vertex in $B$ in $G_1$, and note by construction that $\ell$ will be the degree of a vertex in $C$ in $G_1$.  We now show how to extend each of these eigenvectors to an eigenvector for $G_1$.  We have the following:
\begin{itemize}
\item For an eigenvector of the form ${\bf x}={\bf b}$ (i.e., only nonzero on the vertices of $B$) and $\lambda=0$ then ${\bf y}={\bf b}$ is an eigenvector for $A+tD$ on $G_1$ associated with the eigenvalue of $\gamma=ts$.
\item For an eigenvector of the form ${\bf x}={\bf c}$ (i.e., only nonzero on the vertices of $C$) and $\lambda=0$ then ${\bf y}={\bf c}$ is an eigenvector for $A+tD$ on $G_1$ associated with the eigenvalue of $\gamma=t\ell$.
\item For an eigenvector ${\bf x}={\bf b}+{\bf c}$ associated with an eigenvalue $\lambda\ne0$, we can construct two new eigenvectors, namely
\begin{align*}
&{\bf y}_1=2\lambda{\bf b}+\big(-t(s-\ell)+\sqrt{t^2(s-\ell)^2+4\lambda^2}\big){\bf c}
\mbox{~~and~~}\\
&{\bf y}_2=2\lambda{\bf b}+\big(-t(s-\ell)-\sqrt{t^2(s-\ell)^2+4\lambda^2}\big){\bf c}.
\end{align*}
These are associated with the eigenvalues of $A+tD$ for $G_1$
\[
\gamma_1={t(s+\ell)+\sqrt{t^2(s-\ell)^2+4\lambda^2}\over 2}\quad\mbox{and}\quad
\gamma_2={t(s+\ell)-\sqrt{t^2(s-\ell)^2+4\lambda^2}\over 2}
\]
respectively.

(Note, as in the proof of Theorem~\ref{thm:biswitch}, that we are not creating too many new eigenvectors and eigenvalues, i.e., by Lemma~\ref{lem:bihelp}(b) we can pair up the eigenvalues and eigenvectors.  In particular, ${\bf x}={\bf b}+{\bf c}$ and $\lambda$ will create the same new eigenvectors and eigenvalues as ${\bf x}={\bf b}-{\bf c}$ and $-\lambda$.)
\end{itemize}

The verification that each of these is an eigenvector/eigenvalue pair for $G_1$ reduces to verifying \eqref{eq:t} for each vertex, similarly as was done in Theorem~\ref{thm:biswitch}.  We will skip these routine computations.

Next we note by Lemma~\ref{lem:bihelp}(c) that ${\bf 1}_B$ and ${\bf 1}_C$ are orthogonal to each of the new eigenvectors  for $1<i<m$.  Since the dimension of the space orthogonal to the new eigenvectors restricted to $P_1$ is two, any vector orthogonal to all of these new vectors {\em must}\/ be a linear combination of ${\bf 1}_B$ and ${\bf 1}_C$ when restricted to $P_1$.  In particular, if we let ${\bf y}$ be an eigenvector of $G_1$ that is orthogonal to all of the new eigenvectors, then ${\bf y}$ restricted to $B$ and $C$ must be of the form ${\bf y}|_{B\cup C}=b{\bf 1}_B+c{\bf 1}_C$ for some constants $b$ and $c$.  In other words, the  eigenvector ${\bf y}$ is constant on the vertices of $B$ and constant on the vertices of $C$.

Everything that we have done for $G_1$ carries over to $G_2$.  In particular, since $P_1$ and $P_2$ are cospectral, then the newly found $\gamma_i$ are the same for both graphs, and further, every other eigenvector orthogonal to the ones found must be constant on the vertices of $B$ and the vertices of $C$.  To finish off the first case of the theorem, we now only need to observe that {\em any}\/ eigenvector which is orthogonal to all of the new eigenvectors in $G_1$ is also an eigenvector for $G_2$ for the same eigenvalue.  This again is done in the exact same way as in Theorem~\ref{thm:biswitch} and we will skip the computations here.

In summary, we were able to find $m-2$ eigenvectors of $A+tD$ for $G_1$ and $G_2$ that gave the same set of eigenvalues.  For any other eigenvector orthogonal to these, the same eigenvector worked for both graphs and so the remaining set of eigenvalues also agreed.  So we can conclude that $G_1$ and $G_2$ are cospectral with respect to the matrix $A+tD$.

\bigskip

Now we turn to the case when the vertices on $B$ induce a complete graph.

Using the same notation as in the previous case, we will extend the eigenvectors of $P_1$ for the adjacency matrix to eigenvectors of $H_1$ for $A+tD$.  We have the following:
\begin{itemize}
\item For an eigenvector of the form ${\bf x}={\bf b}$ (i.e., only nonzero on the vertices of $B$) and $\lambda=0$ then ${\bf y}={\bf b}$ is an eigenvector for $A+tD$ on $H_1$ associated with the eigenvalue of $\gamma=ts-1$.
\item For an eigenvector of the form ${\bf x}={\bf c}$ (i.e., only nonzero on the vertices of $C$) and $\lambda=0$ then ${\bf y}={\bf c}$ is an eigenvector for $A+tD$ on $H_1$ associated with the eigenvalue of $\gamma=t\ell$.
\item For an eigenvector ${\bf x}={\bf b}+{\bf c}$ associated with an eigenvalue $\lambda\ne0$, we can construct two new eigenvectors, namely
\begin{align*}
{\bf y}_1&=2\lambda{\bf b}+\big(1-t(s-\ell)+\sqrt{\big(1-t(s-\ell)\big)^2+4\lambda^2}\big){\bf c}\mbox{~~and}\\
{\bf y}_2&=2\lambda{\bf b}+\big(1-t(s-\ell)-\sqrt{\big(1-t(s-\ell)\big)^2+4\lambda^2}\big){\bf c}.
\end{align*}
These are associated with the eigenvalues of $A+tD$ for $H_1$
\begin{align*}
\gamma_1&={t(s+\ell)-1+\sqrt{\big(1-t(s-\ell)\big)^2+4\lambda^2}\over 2}\mbox{~~and}\\
\gamma_2&={t(s+\ell)-1-\sqrt{\big(1-t(s-\ell)\big)^2+4\lambda^2}\over 2}
\end{align*}
respectively.

(Note, as in the proof of Theorem~\ref{thm:biswitch}, that we are not creating too many new eigenvectors and eigenvalues, i.e., by Lemma~\ref{lem:bihelp}(b) we can pair up the eigenvalues and eigenvectors.  In particular, ${\bf x}={\bf b}+{\bf c}$ and $\lambda$ will create the same new eigenvectors and eigenvalues as ${\bf x}={\bf b}-{\bf c}$ and $-\lambda$.)
\end{itemize}

The verification that each of these is an eigenvector/eigenvalue pair for $H_1$ reduces to verifying \eqref{eq:t} for each vertex, similarly as was done in Theorem~\ref{thm:biswitch}.  We will skip these routine computations.

The remainder of this case now proceeds as before.  Namely, everything that we did for $H_1$ carries over for $H_2$ (counting multiplicity of eigenvalues).  Further, any other eigenvector orthogonal to the ones given must be constant on $B$ and $C$, so any other eigenvector for $A+tD$ which works for $H_1$ also works for $H_2$.  Therefore we can conclude that they have the same set of eigenvalues, i.e., are cospectral with respect to $A+tD$.
\end{proof}

\section{Large families of mutually cospectral graphs}\label{sec:applications}
We can use Theorem~\ref{thm:biswitch} to construct large families of mutually cospectral graphs for the normalized Laplacian that need not be bipartite nor have the same degree sequences. 

\begin{example}
Let $n$, $k$ be positive integers such that $m_1+\cdots+m_k=n$ is a partition of $n$ into $k$ positive integer parts.  Let $FB(m_1,m_2,\ldots,m_k)$ be the graph on $n+k$ vertices $b_1,\ldots,b_n,v_1,\ldots,v_k$, where the $b_i$ induce a complete graph, each $v_i$ is only adjacent to exactly $m_i$ vertices from $b_1,\ldots,b_n$, and each $b_i$ is adjacent to exactly one of the $v_i$.  The collection of all $FB(m_1,m_2,\ldots,m_k)$ for partitions $m_1+m_2+\cdots+m_k=n$ give a family of cospectral graphs for the normalized Laplacian, i.e., the spectrum is completely determined for the normalized Laplacian by $n$ and $k$.
\end{example}

\begin{figure}[hftb]
\centering
\includegraphics[scale=0.9]{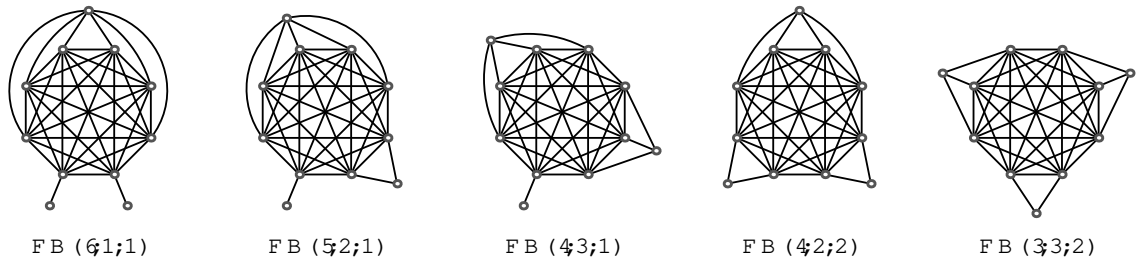}
\caption{Fuzzy ball graphs for $n=8$ and $k=3$.}
\label{fig:fuzzy}
\end{figure}

We have dubbed these graphs ``fuzzy balls'' since they consist of a complete graph (the $b_i$) with some extra vertices attached (the fuzz).  In Figure~\ref{fig:fuzzy}, we have shown the cospectral family that corresponds to $n=8$ and $k=3$.

To see that these are all cospectral, we observe by Theorem~\ref{thm:biswitch} that $FB(m_1,m_2,\ldots,m_k)$ is cospectral with $FB(m_1+m_i-1,m_2,\ldots,m_{i-1},1,m_{i+1},\ldots,m_k)$ by using the pair of cospectral graphs $K_{1,m_1}\cup K_{1,m_i}$ and $K_{1,m_1+m_i-1}\cup K_{1,1}$ where $B$ induces a complete graph.  Therefore, applying this idea $k-1$ times, we can conclude that for each partition of $n$ into $k$ parts we have $FB(m_1,m_2,\ldots,m_k)$ is cospectral with $FB(n-k+1,1,\ldots,1)$.  In particular, they are all mutually cospectral.

\begin{example}
Let $n$, $k$ be positive integers such that $m_1+\cdots+m_k=n$ is a partition of $n$ into $k$ positive integer parts.  Let $IS(m_1,m_2,\ldots,m_k)$ be the graph on $n+k+1$ vertices $a$, $b_1$, $\ldots$, $b_n$, $v_1$, $\ldots$, $v_k$, where the $b_i$ and $a$ induce a star graph with $a$ as the central vertex, each $v_i$ is only adjacent to exactly $m_i$ vertices from $b_1,\ldots,b_n$, and each $b_i$ is adjacent to exactly one of the $v_i$.  The collection of all $IS(m_1,m_2,\ldots,m_k)$ for partitions $m_1+m_2+\cdots+m_k=n$ give a family of cospectral graphs, i.e., the spectrum is completely determined by $n$ and $k$.
\end{example}

\begin{figure}[hftb]
\centering
\includegraphics[scale=0.9]{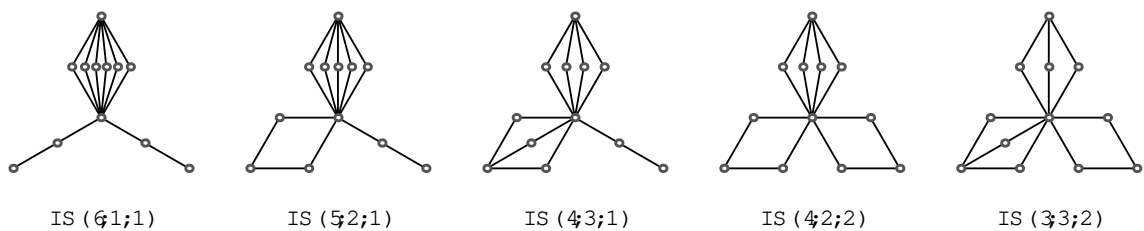}
\caption{Inflated star graphs for $n=8$ and $k=3$.}
\label{fig:stars}
\end{figure}

We have dubbed these graphs ``inflated stars''.  In Figure~\ref{fig:stars}, we have shown the cospectral family that corresponds to $n=8$ and $k=3$.

Again, to see that these are all cospectral, we observe by Theorem~\ref{thm:biswitch} that the graph $IS(m_1,m_2,\ldots,m_k)$ is cospectral with $IS(m_1+m_i-1,m_2,\ldots,m_{i-1},1,m_{i+1},\ldots,m_k)$.  Therefore, applying this idea $k-1$ times, we can conclude that for each partition of $n$ into $k$ parts, we have $IS(m_1,m_2,\ldots,m_k)$ is cospectral with $IS(n-k+1,1,\ldots,1)$.  In particular, they are all mutually cospectral.

In both of these examples, we are using partitions to form large families.  The number of partitions grows subexponentially with $n$, however it is possible to use Theorems~\ref{thm:biswitch} or~\ref{thm:biregswitch} to construct large mutually cospectral families that grow exponentially with $n$, either with respect to the normalized Laplacian, or more generally with respect to $A+tD$.

\begin{theorem}\label{thm:nlfam}
For $n$ large, there exists a family of $2^{\lfloor n/7\rfloor}$ non-isomorphic, mutually cospectral graphs with respect to the normalized Laplacian.
\end{theorem}
\begin{proof}
Find a graph $G$ on $\lfloor n/7\rfloor$ vertices that has a trivial automorphism group (for $n\ge 42$ this is easy to do).  For each vertex $v$ of $G$, attach one of the two inflated stars shown in Figure~\ref{fig:nlfam} by identifying $v$ with the vertex marked $v$ in the inflated star.  This constructs $2^{|G|}=2^{\lfloor n/7\rfloor}$ graphs on $7\lfloor n/7 \rfloor\le n$ vertices.  
\begin{figure}[ht]
\centering
\includegraphics{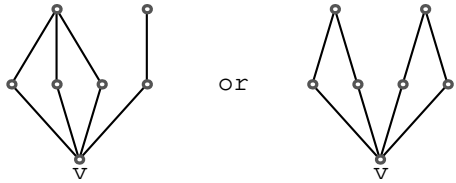}
\caption{The two different inflated stars to add to vertices $v$ in Theorem~\ref{thm:nlfam}.}
\label{fig:nlfam}
\end{figure}
We note that for any one of these graphs we can easily recover $G$ and identify which attachment was done at $v$ and so these graphs are all non-isomorphic.  Further, by Theorem~\ref{thm:biswitch} changing our choice of widget at any single vertex does not change the spectrum and so all of the graphs are cospectral.
\end{proof}

We have not tried to optimize the construction.  For example, if instead of using the family of inflated stars with $n=4$ and $k=2$ we used the family of inflated stars with $n=21$ and $k=5$ (of which there are $101$ such non-isomorphic members), then we get a family with rate of growth $101^{\lfloor n/26\rfloor}\approx 1.19423^n$ which is faster than $2^{\lfloor n/7\rfloor}\approx 1.10409^n$.

We similarly have the following result for graphs cospectral with respect to $A+tD$.

\begin{theorem}\label{thm:atdfam}
For $n$ large, there exists a family of $2^{\lfloor n/13\rfloor}$ non-isomorphic, mutually cospectral graphs with respect to the matrix $A+tD$.
\end{theorem}
\begin{proof}
Find a graph $G$ on $\lfloor n/13\rfloor$ vertices that has a trivial automorphism group (for $n\ge 78$ this is easy to do).  For each vertex $v$ of $G$ attach one of the two graphs shown in Figure~\ref{fig:atdfam}, which are taken from Figure~\ref{fig:2bi}.  This constructs $2^{|G|}=2^{\lfloor n/13\rfloor}$ graphs on $13\lfloor n/13 \rfloor\le n$ vertices.  
\begin{figure}[ht]
\centering
\includegraphics{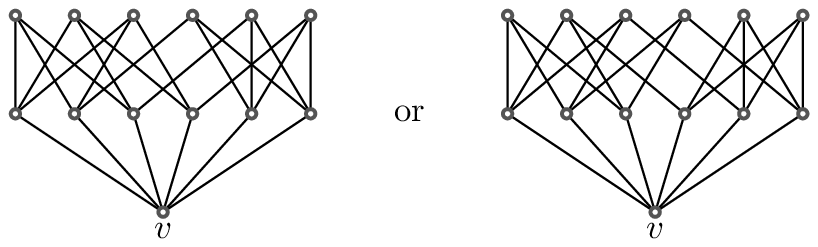}
\caption{The two different widgets to add to vertices $v$ in Theorem~\ref{thm:atdfam}.}
\label{fig:atdfam}
\end{figure}
We note that for any one of these graphs we can easily recover $G$ and identify which attachment was done at $v$, so these graphs are all non-isomorphic.  Further, by Theorem~\ref{thm:biregswitch} changing our choice of widget at any single vertex does not change the spectrum and so all of the graphs are cospectral.
\end{proof}

\section{Concluding remarks}\label{sec:comments}
The method of swapping bipartite subgraphs was discovered by examining all cospectral pairs of graphs on at most eight vertices (all cospectral pairs on at most seven vertices has previously appeared in Tan \cite{Tan}), which included several simple examples of this type.  The set of all cospectral graphs on at most eight vertices can be easily constructed in {\sc sage} (see Appendix~\ref{sagecode}).

\begin{figure}[ht]
\centering
\includegraphics{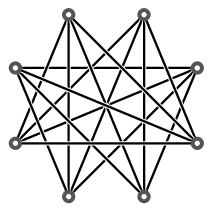}\hspace{.5in}\includegraphics{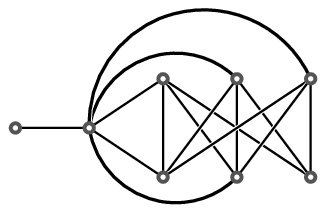}
\caption{A regular non-bipartite graph which is cospectral for the normalized Laplacian with a graph which is not regular.}
\label{fig:regular}
\end{figure}

By looking at the small cases, we also discovered some new graphs which are regular, but are cospectral for the normalized Laplacian with graphs which are not regular.  We have already seen that the four-cycle and the eight-cycle are members of  cospectral pairs (see Figures~\ref{fig:1d} and \ref{fig:other1}).  Also, $K_{n,n}$ is cospectral with $K_{p,q}$ for $p+q=2n$.  Another example of a graph which is regular, but is cospectral with a graph which is not regular, is shown in Figure~\ref{fig:regular}.  These examples show that the normalized Laplacian cannot, in general, detect whether a graph is regular.  

\begin{figure}[ht]
\centering
\hfil \subfloat[Cospectral pair]{\includegraphics{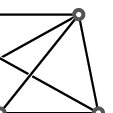}\hspace{.5in}\includegraphics{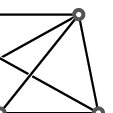}\label{fig:4a}}
\hfil \hfil \subfloat[Cospectral pair]{\includegraphics[angle=90]{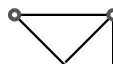}\hspace{.5in}\includegraphics[angle=90]{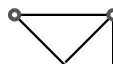}\label{fig:4b}}
\hfil

\hfil \subfloat[Cospectral pair]{\includegraphics{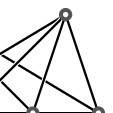}\hspace{.5in}\includegraphics{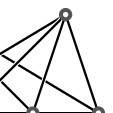}\label{fig:4c}}
\hfil \hfil \subfloat[Cospectral pair]{\includegraphics[angle=90]{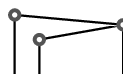}\hspace{.5in}\includegraphics[angle=90]{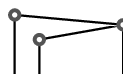}\label{fig:4d}}
\hfil

\caption{Examples of graphs which are cospectral with a subgraph for the normalized Laplacian.}
\label{fig:sub}
\end{figure}

\begin{figure}[ht]
\centering
\hfil \subfloat[Cospectral pair]{\includegraphics{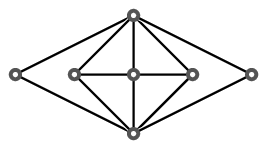}\hspace{.5in}\includegraphics{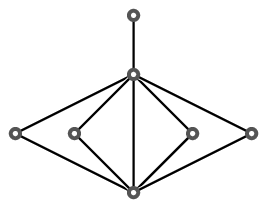}\label{fig:5a}}
\hfil \hfil \subfloat[Cospectral pair]{\includegraphics{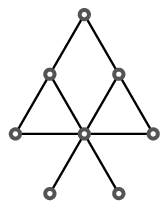}\hspace{.5in}\includegraphics{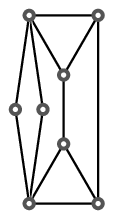}\label{fig:5b}}
\hfil

\hfil \subfloat[Cospectral pair]{\includegraphics{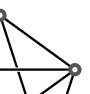}\hspace{.5in}\includegraphics{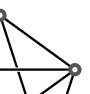}\label{fig:5c}}
\hfil \hfil \subfloat[Cospectral pair]{\includegraphics{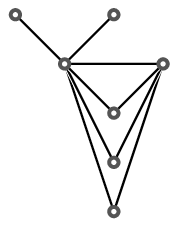}\hspace{.5in}\includegraphics{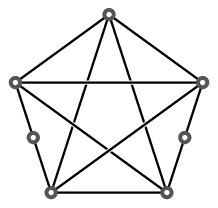}\label{fig:5d}}
\hfil

\caption{Examples of graphs which are cospectral for the normalized Laplacian but differ in the number of edges.}
\label{fig:dif}
\end{figure}

As noted in the introduction, the number of edges in a graph is not uniquely determined by the spectrum of the normalized Laplacian.  As a consequence, it is possible for a graph to be cospectral with one of its subgraphs; examples of this are shown in Figure~\ref{fig:sub}.  Further examples of graphs which are cospectral with respect to the normalized Laplacian but have differing number of edges are shown in Figure~\ref{fig:dif} (there is currently no known method to generate cospectral graphs with differing number of edges other than using complete bipartite graphs).  This also indicates the difficulty in counting the number of cospectral graphs.  While for the adjacency and the combinatorial Laplacian, we could first subdivide the graphs according to the number of edges and work with this coarsening when finding cospectral pairs (see \cite{HaemersSpence}), this is no longer possible for the normalized Laplacian.

It would be interesting to find techniques for constructing cospectral graphs for the normalized Laplacian which differ in the number of edges.  Similarly, it would be interesting to investigate the cospectral graphs for the signless Laplacian and see if there are any constructions unique to forming cospectral graphs for that matrix (similar to what Theorem~\ref{thm:biswitch} is for the normalized Laplacian).

\appendix

\section{{\sc sage} code to generate cospectral graphs}
\label{sagecode}

The code listing below shows how to use Sage version 4.5.2 \cite{sage} to generate and check cospectral graphs with 8 vertices.  In the listing below, we specifically generate cospectral graphs with respect to the normalized Laplacian.  The lines starting with the number sign ``\#'' are comment lines.

\lstset{language=Python, mathescape=true,numbers=left,
numberstyle=\tiny, stepnumber=5, numberfirstline=false, firstnumber=1,
basicstyle=\small\ttfamily,columns=fixed, commentstyle=\textsl,showstringspaces=false}

\begin{lstlisting}
def DinverseA(g):
    "Calculate D^{-1}A for a graph g"
    A=g.adjacency_matrix().change_ring(QQ)
    for i in range(g.order()):
        A.rescale_row(i, 1/len(A.nonzero_positions_in_row(i)))  
    return A

# Calculate the graphs and store them in cospectral_list.
# We use the DinverseA function defined above to determine
# if two graphs are cospectral.
# This command takes a few minutes to complete.
cospectral_list=graphs.cospectral_graphs(8,
                          graphs=lambda g: min(g.degree())>0)

# Give a list of all graph6 strings of these graphs.
graph6_list=[[g.graph6_string() for g in glist] 
             for glist in cospectral_list]

# Show all of the graphs, each row being cospectral.
for glist in cospectral_list: show(glist)

# Get the first two cospectral graphs and check
# to make sure they are cospectral.
graph1, graph2 = cospectral_list[0]
graph1_poly=graph1.laplacian_matrix(normalized=True).charpoly()
graph2_poly=graph2.laplacian_matrix(normalized=True).charpoly()
graph1_poly==graph2_poly

# Construct fuzzy ball graphs corresponding to the partitions 
# 4=3+1 and 4=2+2 and check that they are cospectral.
fb1=graphs.FuzzyBallGraph([3,1],0)
fb2=graphs.FuzzyBallGraph([2,2],0)
fb1_poly=fb1.laplacian_matrix(normalized=True).charpoly()
fb2_poly=fb2.laplacian_matrix(normalized=True).charpoly()
fb1_poly==fb2_poly
\end{lstlisting} 

\end{document}